\documentclass[11pt]{article}

\usepackage{amsmath, amssymb}
\usepackage{graphicx}
\usepackage{float,here}
\usepackage{xcolor}
\usepackage{dsfont}
\usepackage{mathrsfs}
\usepackage{hyperref}

\hypersetup{
colorlinks=true, 
breaklinks=true, 
urlcolor= blue, 
linkcolor= blue, 
bookmarksopen=true, 
pdftitle={}, 
pdfauthor={Christian Bourdarias}, 
pdfsubject={} 
}

\newcommand{\cqfd}{{\nobreak\hfil\penalty50\hskip2em\hbox{}\nobreak\hfil
$\square$\qquad\parfillskip=0pt\finalhyphendemerits=0\par\medskip}}

\newcommand{\R}{\mathbb{R}}

\newcommand{\N}{\mathbb{N}}

\def\abs#1{\vert#1\vert}

\def\norminf#1{\parallel#1\parallel_{\infty}}

\newcommand{\pr}{{\bf \textit{Proof: }}}

\newcommand{\1}{\mathds{1}}

\newtheorem{theorem}{Theorem}[section]
\newtheorem{corollary}{Corollary}[section]
\newtheorem{proposition}{Proposition}[section]
 \newtheorem{lemma}{Lemma}[section]

\newtheorem{remark}{Remark}[section]

\title{ \bf Kinetic formulation of a $2\times2$ hyperbolic system  arising in  gas chromatography}

\author {C. Bourdarias
\thanks{U.  Savoie Mont Blanc, LAMA, UMR CNRS 5127, 73376 Le Bourget-du-Lac,
bourdarias@univ-savoie.fr},
M. Gisclon
\thanks{U.  Savoie Mont Blanc, LAMA, UMR CNRS 5127, 73376 Le
Bourget-du-Lac, gisclon@univ-savoie.fr}
and S. Junca
\thanks{U. de Nice, Labo. JAD, UMR CNRS 7351, Parc Valrose,
06108 Nice,  junca@unice.fr, and   INRIA  M\'editerran\'ee, Team Coffee, 2004 route des Lucioles, 06902, Sophia-Antipolis, France}
}

\date{\today}

\begin{document}
\everymath{\displaystyle}
\noindent

\bibliographystyle{plain}

\maketitle

\tableofcontents

\abstract{ The PSA system commonly used in the context of gas-solid chromatography is reformulated as a single kinetic equation using an additional kinetic variable. A kinetic numerical scheme  is built from this new formulation and its behavior is tested on solving the Riemann problem in different configurations leading to single or composite waves. }\\[2mm]


\section{Introduction}

Since the work of P.-L. Lions, B. Perthame and E. Tadmor (\cite{LPT94,LPT}), it is well known that 
multidimensional scalar conservation laws and some systems can be formulated as a kinetic equation using an additional kinetic variable. 
The so-called kinetic formulation of nonlinear hyperbolic systems of conservation laws  reduces them to a linear equation  on a nonlinear quantity related to the conservative unknowns, moreover it allows to recover  all the entropy inequalities. It turns out to be a powerful
tool to derive mathematical properties such that regularizing effects or compactness results and also efficient numerical schemes.
The method was used by several authors who gave further examples of kinetic formulations: the  system of chromatography (\cite {JPP95,JPP96}), the Shallow Water system (\cite{PS01}) for instance.
The objective of the present work is to apply the machinery of entropy to derive a kinetic formulation for the the so-called PSA system introduced in Section \ref{psa} and already studied by the authors from various points of view (\cite{BGJ1,BGJ2,BGJ5,BGJP}). As a first application, we construct a kinetic numerical scheme, state some of its properties and test it by solving the Riemann problem.

The paper is organized as follows.
In Section \ref{psa}, we present the PSA  system and in Section \ref{sh}, we recall basics results of hyperbolicity and entropies. These two sections summarize the results of previous work by the same authors essential to understanding the problem addressed and contain no new result.\\
  In Section \ref{kfPSA}, we build and analyze a kinetic formulation of System (\ref{sysad}). In Section \ref{ksPSA}  a kinetic scheme is built using the preceding formulation. The last section is devoted to the numerical validation of the scheme based on the resolution of the Riemann Problem.


\section{The PSA system} \label{psa}

Pressure Swing Adsorption (PSA) is a technology that is used to separate some species from a gas under pressure according to these species'   molecular characteristics and affinity for an adsorbent material. PSA is used extensively in the production and purification of oxygen, nitrogen and hydrogen for industrial uses. It can also be used to separate a single gas from a mixture of gases. A typical PSA system involves a cyclic process where a number of connected vessels containing adsorbent material undergo successive pressurization and depressurization steps in order to produce a continuous stream of purified product gas.

As in previous papers by the authors on this subject, we focus on a model describing a step of the cyclic process, restricted to isothermal behavior. As in general fixed bed chromatography, each of the $d$ species ($d\geq 2$) simultaneously exists under two phases, a gaseous and movable one with concentration $c_i (t,x)$ or a solid (adsorbed) other with concentration $q_i (t,x)$, $1\leq i\leq d$. Moreover it is assumed that these concentations are at equilibrium, i.e. $q_i(t,x)=q_i^*(c_1,c_2)$, where the so-called isotherms $q_i^*$ satisfy  
\begin{equation}\label{isoth}
\frac{\partial q_i^*}{\partial c_i} \geq 0, \quad i \in \{1,2\}.
\end{equation}
In gas chromatography, velocity variations accompany changes in gas composition, especially in the case of high concentration solute: it is known as the sorption effect. This effect is taken into account through a constraint on the pressure, assumed to be constant. The reader can refer for instance to \cite{DCRBT88}: ``Fixed-Bed Adsorption of Gases : Effect of Velocity Variations on Transition Types". 

In \cite{BGJP},  the original model, which is  nothing else that material balances for two adsorbable components, is written into the dimensionless form:

\begin{eqnarray}
        \partial_t ( c_1 + q_1^*(c_1,c_2)  ) +   \partial_x  ( c_1 \, u )   & = &  0,  \label{ad1}\\
          \partial_t ( c_2 + q_2^*(c_1,c_2)  ) +   \partial_x  ( c_2 \, u )   & = &  0, \label{ad2}\\
          c_1 + c_2    &  = &  1.\label{ad3} 
\end{eqnarray}

In this isothermal model, the constraint of constant pressure is achieved through Eq. \eqref{ad3}. In the case where  an instantaneous equilibrium is not assumed, the corresponding system was studied  from both theoretical and numerical points of view by  Bourdarias \cite{B92,B98}.

 Setting $c=c_1\in[0,1]$ (then $c_2=1-c$), $q_i(c)=q_i^*(c,1-c)$, $i=1,2$, and adding Eqs. \eqref{ad1}-\eqref{ad2}, we obtain finally the following system for $x>0$ and $t>0$ which is written in a form ($x$-derivative first)  justified in the next section:
\begin{equation}\label{sysad}
\left\{   \begin{array}{ccl}
     \partial_x  (u\, c)+  \partial_t I(c)& = & 0, \\
   \partial_x u+  \partial_t h(c)  & = & 0, 
  \end{array}\right.
\end{equation}
where
$$    h(c)=q_1(c)+q_2(c) \geq 0, \quad I(c)= c + q_1(c).  $$
With theses notations,  the relations \eqref{isoth} read $q_1'\geq 0 \geq q_2'$.

Following \cite{DCRBT88}, we introduce a  function which will plays a central role in the nonlinear study of the system, namely,
\begin{equation}\label{f}
  f(c) = c_2 \, q_1 - c_1 \, q_2  =      q_1 (c) - c \, h(c).  
\end{equation}

We will also make use of following functions only depending on the
isotherms \cite{BGJ2}: 
\begin{itemize}
\item
$H(c)=  1 + q_1' -c h' =1+(1-c)q_1'-cq_2' \geq 1$, 
\item
$G(c) = \exp g(c)$ where   $g' = -\displaystyle  \frac{h'}{H}$.
\end{itemize}


\section{Hyperbolicity and entropies }\label{sh}


For self contain, we recall without proofs some results exposed in \cite{BGJ2}.\\
As pointed out by Rouchon and \textit{al.} (\cite{RSVG88}), it is possible to analyze System
(\ref{sysad}) in terms of {\bf hyperbolic} system of P.D.E. provided the time and space
variables are exchanged : this is why System \eqref{sysad} is presented under this unususal form.
 In this framework, the vector state is  
$U=\left(
\begin{array}{l}
u \\
m
\end{array}\right)$ 
where $m=u\,c $ is the flow rate of the first species. In this vector, $u$ must be understood as
$u\,(c_1+c_2)$, that is the total flow rate. The initial-boundary value  problem is
then System \eqref{sysad} for $x>0$ and $t>0$ supplemented by the initial ($x=0$) and boundary data ($t=0$): 
\begin{equation} \label{sysad0}
\left\{   \begin{array}{ccl}
 \vspace{2mm}c(t,0) &=&c_b(t)   \in [0,1],\quad t>0,\\
 u(t,0)&=&u_b(t)>0,\quad t>0\\\\
 \vspace{2mm}c(0,x)&=&c_0(x) \in [0,1],  \quad x > 0,
 \end{array}\right.
\end{equation}
where   
\begin{eqnarray*}
    & 0 \leq  c_b, c_0  \leq 1,    \label{0c1}& \\
    &0 <  \inf_\R u_b  \leq \sup_\R u_b  < + \infty . & \label{infsupub}
\end{eqnarray*}
For this system,  the first two equations of \eqref{sysad0} correspond to the initial data and the last one to the boundary data. That is to say that the variable $x$ is progressive, i.e. time-like, and  $t$ is a space-like variable. To be clear,  we distinguish the physical time $t$ to the mathematical time or hyperbolic time $x$. The mathematical initial value problem is physically relevant for applications because experimenters only control $c_b, u_b$, and $c_0$ can be viewed as an equilibrium reached before the beginning of the process.\\
The eigenvalues of the Jacobian matrix of the flux are $0$ and $\lambda=\frac{H(c)}{u}$, thus the system is strictly hyperbolic as long as $u>0$: we will show in the last section that it is ensured for the solution of the Riemann Problem thanks to the assumption $\inf_\R u_b>0$. 
Moreover $\lambda$ is genuinely nonlinear in each domain where $f"\neq 0$. The Riemann invariants are $c$  and $W=u\,G(c)$ associated to the eigenvalues $0$ and $\lambda$ respectively.\\
We have shown in \cite{BGJ2} that there are two families of entropies:
$u\,\psi(c)$ and $\phi(u\,G(c))$, where $\phi$ and $\psi$ are any real smooth functions. The corresponding entropy flux $Q(c)$ of the first family satisfies
\begin{eqnarray*}
Q'(c)  & = & h'(c)\, \psi(c) +  H(c)\, \psi'(c) .
\end{eqnarray*}
The first family is degenerate convex (in variables  $(u,uc)$) provided $\psi''\geq 0$. 
So we seek entropy solutions which satisfy 
$$  \partial_x \left(u\,\psi(c) \right) + \partial_t Q(c) \leq 0,  $$
in the distribution sense. The second family is not always convex. There are only two
interesting cases where this family is convex, namely $\pm G''(c) > 0$ for all $c \in [0,1]$.
When $G''>0$ and $\alpha > 1$, we expect to have 
 $\partial_x ( u\,G(c))^\alpha \leq 0 $ which reduces to 
  $\partial_x ( u\,G(c)) \leq 0 $.
In the same way, if $ G''<0$, we get $ \partial_x ( u\,G(c)) \geq 0 $.


\section{Kinetic Formulations of PSA System}\label{kfPSA}


In this section, we consider weak solutions of System \eqref{sysad} and we give two kinetic
formulation. This requires the knowledge of a complete family of supplementary conservation
laws or more precisely the weak entropy inequalities (\cite{PT1}).   The first general formulation will be used later to build a kinetic scheme. The second formulation restricted to  a convex assumption on isotherms is just mentionned in the last subsection.


\subsection{Main kinetic inequality}

To have a general kinetic formulation we use the family $u\,\psi(c)$, where $\psi'' \geq 0$. Despite the fact that this family is
always degenerate convex, this family has the great advantage to be convex without convex assumption
on isotherms.

Let us introduce the  classic function $\chi$:

\begin{eqnarray}  \label{kai}
 \chi(c,\xi) =\1_{]0,c[}(\xi)=\left \{ \begin{array}{cl} 
                1  & \mbox{ if } 0 < \xi < c, \\
               0  & \mbox{else}.
             \end{array} \right.
\end{eqnarray}

This function enjoys the following simple properties:
$$\mbox{ supp}(\chi)=[0,c],\quad \int_\R \chi(c, \xi) \, d \xi=c$$ and $$\forall g \in {\cal C}^1(\R), \
\int_\R g'(\xi) \, \chi(c,\xi) \, d \xi =g(c)-g(0).$$
 Since $ c \in (0,1)$ we define $\chi$ only for $(c,\xi) \in (0,1)\times (0,1)$ and we have also $ \chi(c,\xi) =\1_{]\xi,1[}(c) $. 
  
Moreover, this function satisfies the fundamental Gibbs property (also called Brenier's Lemma):
let $S$ be a convex function, $c \in \R$ and consider the minimization problem

\begin{equation} \label{kf}
\inf\{ \int_0^1 S'(\xi)\,\phi(c,\xi) \, d\xi; \ \phi\in L^1_\xi(0,1), \  0 \leq \phi\leq
1
         \mbox{ and }  \int_0^1 \phi(c,\xi) \, d\xi = c    \}.
\end{equation} 

Y. Brenier (\cite{P1,Br1})) showed that the minimization problem achieves its
minimum at $f=\chi$ and that if $S$ is strictly convex the minimizer is unique.

\bigskip

In order to derive the kinetic formulation of  (\ref{sysad}), we need to state two technical but easy results.

\begin{lemma}\label{lem1}
The distribution $\partial_c \chi$ satisfies: $\forall \phi \in C^\infty_0((0,1)^2)$,
$<\phi, \partial_c \chi> = \int_0^1 \phi (\xi,\xi) \, d \xi$.
\end{lemma}

\pr for all $\phi \in C^\infty_0((0,1)^2)$ we have, using $ \chi(c,\xi) =\1_{]\xi,1[}(c) $ :

\begin{eqnarray*}
<\phi, \partial_c \chi>=-< \partial_c\phi,\chi>&=&- \int_0^1 \int_0^1 \partial_c \phi(c, \xi)\,\chi(c,\xi) \, d \xi \, dc \\
&=& - \int_0^1 \left( \int_\xi^1 \partial_c \phi(c, \xi)\, dc \, \right) d \xi\\
&=&\int_0^1 \phi(\xi,\xi) \, d \xi.
\end{eqnarray*}
\cqfd
In the following  lemma, $H$ is the function defined at the end of the second section.

\begin{lemma}\label{lem2}
The function $P(c,\xi)=H(\xi)\, \chi(c,\xi)$ satisfies $\partial_c P= H(c) \,\partial_c \chi$.
\end{lemma}

\pr
 on one hand, using Lemma \ref{lem1},  
$$<\phi(c,\xi), H(c) \, \partial_c \chi>=
<H(c) \, \phi(c,\xi), \partial_c \chi >=\int_0^1 H(\xi) \, \phi(\xi,\xi) \, d
\xi.$$
On the other hand, 
\begin{eqnarray*}
<\phi(c,\xi), \partial_c P> =-<\partial_c \phi(c,\xi), P>&=& -<\partial_c \phi(c,\xi), H(\xi) \, \chi(c,\xi)>\\
&=&-\int_0^1 \int_0^1 H(\xi) \, \partial_c \phi(c,\xi) \, \chi(c,\xi) \, dc \, d \xi\\
&=& -\int_0^1 H(\xi) \left( \int_\xi^1 \partial_c \phi(c, \xi) \, dc \right) \, d \xi \\
&=&\int_0^1 H(\xi) \, \phi(\xi,\xi) \, d\xi.
\end{eqnarray*}
\cqfd 

We are now ready to give our main result:

\begin{theorem}\label{kinform}
If $(u,c)$ is a weak entropy solution of System (\ref{sysad}), then there exists a
non\-ne\-ga\-ti\-ve measure $m(t,x,\xi)$ such that: 
\begin{equation} \label{eqkf}
\partial_x( u \, \chi(c,\xi)) + a(\xi) \, \partial_t \, \chi(c,\xi) 
+ \partial_t \left(h(c) \, \chi(c,\xi)  \right)= \partial_\xi \, m,
\end{equation}
where $a$ is given by $a(\xi)=H(\xi)-h(\xi)= 1+f'(\xi)$.
\end{theorem}

\pr
we begin to obtain the kinetic formulation (\ref{eqkf}) by writing entropy inequalities for all
$\psi$ such that $\psi'' \geq 0$:

\begin{equation*} \label{iepsi}
\partial_x \left(  u \, \psi(c)\right) + \partial_t Q(c) \leq 0, 
\quad \mbox{ where }  Q'(c)=H(c) \, \psi'(c) + h'(c) \, \psi(c). 
\end{equation*}

With Kruzkhov entropies $\psi(c,\xi)=|c-\xi| -|\xi|$  which satisfy 
$$-\frac{1}{2}\partial_\xi \psi(c,\xi) = \chi(c,\xi),$$
we have a nonnegative measure  $m$ such that:

\begin{equation} \label{KE}
\partial_x \left( u \, \psi(c,\xi)\right) + \partial_t Q(c,\xi) =-2 \, m(t,x,\xi). 
\end{equation}

Applying  $-\frac{1}{2} \partial_\xi$ on the previous equation we get:

\begin{equation*}
\partial_x \left( u \, \chi(c,\xi)\right) + \partial_t \left( -\frac{1}{2} \, \partial_\xi
Q(c,\xi)\right) =\partial_\xi \, m(t,x,\xi). 
\end{equation*} 

So, we have to compute $ -\frac{1}{2} \, \partial_\xi \, Q(c,\xi) $.
Since $$\partial_c \, Q(c,\xi)=  H(c) \, \partial_c \psi(c,\xi) + h'(c)\,\psi(c,\xi),$$
we have, applying  $-\frac{1}{2} \partial_\xi$ and using Lemma \ref{lem2}:

\begin{eqnarray*}
-\frac{1}{2} \, \partial_\xi \, \partial_c Q(c,\xi) &= &H(c) \, \partial_c \chi(c,\xi) + h'(c) \, \chi(c,\xi)\\
&=& \partial_c (H(\xi)\, \chi(c,\xi))+ h'(c) \, \chi(c,\xi).
\end{eqnarray*}

Notice that 
$$\int_0^c h'(y) \, \chi(y,\xi) \, dy = (h(c)-h(\xi)) \, \chi(c,\xi),$$
thus we have:

\begin{equation*}
-\frac{1}{2}\partial_\xi  Q(c,\xi) = 
-\frac{1}{2}\partial_\xi Q(0,\xi) + H(\xi) \, \chi(c,\xi)+ (h(c)-h(\xi)) \,\chi(c,\xi).
\end{equation*}

We get finally

\begin{equation*}
\partial_x( u  \, \chi(c,\xi)) + H(\xi) \, \partial_t \chi(c,\xi) 
+ \partial_t [(h(c) -h(\xi)) \, \chi(c,\xi)]= \partial_\xi m.
\end{equation*}

This is valid with the $\chi$ function associated to the special entropy $\psi(c,\xi)=|c-\xi| -|\xi|$.  
Since Kruzkhov entropies generate all convex functions by convex combinations and density arguments, Theorem  \ref{kinform} holds.
\cqfd

\bigskip

The conversely of previous theorem is the following one:

\begin{theorem}
If there exist  a positive function $u$ such that $\ln u \in L^\infty$, a $\chi$-function  $\phi(t,x,\xi)=\chi(c(t,x),\xi)$ for some function $c$
 and a nonnegative measure  $m$ such that
 \begin{equation}\label{recip}
\partial_x  (u\, \phi(t,x,\xi)) + a(\xi)\,\partial_t\phi(t,x,\xi)+ \partial_t \left(h(c) \, \phi(t,x,\xi)  \right)= \partial_\xi m, 
 \end{equation}
then $(u,c) $ is a weak entropy solution of System (\ref{sysad}).
\end{theorem}

\pr
multiplying Equality (\ref{recip}) by $\psi'$  and integrating over $(0,\xi)$ we get (\ref{KE}). 
With (\ref{KE}) we recover easily System (\ref{sysad}):
first, using  $\psi\equiv \pm 1$ then $Q'(c)=\pm h'(c)$ and we recover the second equation 
of  (\ref{sysad}), next the  choice $\psi\equiv \pm c$ gives $Q'(c)=\pm I'(c)$ and we recover the first
equation of  (\ref{sysad}).
\cqfd

\begin{remark}\rm
Writing Eq. (\ref{eqkf}) under the form
$$\partial_x( u \, \chi(c,\xi)) + \, \partial_t \,\left( (a(\xi) + h(c) ) \, \chi(c,\xi)  \right)=
\partial_\xi \, m$$
we highlight an advection velocity $a(\xi) + h(c)$ which is is not purely kinetic, as in \cite{LPT}.
\end{remark}

For weak solution we can bound the measure $m$ with respect to $L^\infty$ bound of $u$.

\begin{proposition}[A priori bound for defect measure]~\\
If $(u,c)$ is a weak entropic solution of System (\ref{sysad}) and let  $m(t,x,\xi)$ the defect
measure satisfying  kinetic formulation (\ref{eqkf}), then there exists a constant $\alpha > 0$ 
depending on  $||h||_{\infty}$  such that:
$$
\forall T>0,\, \forall X>0,\quad    \int_0^T\int_0^X\int_0^1 m(t,x,\xi) \, dt \, dx \, d\xi 
  \leq  \alpha X + T \| u_b \|_{L^\infty(0,T)}. 
$$
\end{proposition}

\pr
multiplying (\ref{eqkf}) by $S'(\xi)$ such that $S(0)=0$,  we get

\begin{equation*}
\partial_x( u \, S'(\xi) \, \chi(c,\xi))+ \partial_t \left( S'(\xi) \, a(\xi) \,  \chi(c,\xi)
\right)+ \partial_t \left(h(c) \, S'(\xi) \, \chi(c,\xi)  \right)=S'(\xi)  \,\partial_\xi \,
m(t,x,\xi).
\end{equation*}

Integrating by parts the previous equality over $(0,T)_t\times(0,X)_x\times (0,1)_\xi $ we have:

\begin{eqnarray*} \label{eqentfk}
  &  -\int_0^T\int_0^X \int_0^1 S''(c) \ m \ dt \ dx \ d\xi  
\\ 
 = & 
\int_0^T [u(t,X)S(c(t,X)) -u_b(t)S(c_b(t))]dt 
+ \int_0^X \left( A(c)(T,x) -A(c_0)(x)\right)dx 
\end{eqnarray*}

where $A(c)= \int_0^c a(\xi) \, S'(\xi)\, d\xi + h(c) \, S(c).$ Since $ 0\leq c \leq 1$, we have
only to estimate $u(t,X)$. To control $ \int_0^T u(t,X) \,dt$ we use second equation of System
(\ref{sysad}): $\partial_x u = -\partial_t h(c)$ and nonnegativity of $h$, then:

\begin{eqnarray*}
 \int_0^Tu(t,X) \, dt & = & 
    \int_0^Tu_b(t) \, dt + \int_0^X( h(c_0(x))-h(c(T,x))) \, dx  \\
   &\leq &
    \int_0^Tu_b(t) \, dt + \int_0^X h(c_0(x)) \, dx. 
\end{eqnarray*}

With $S(c) =c^2/2$, a constant $\alpha$ depending only on the supremum of $A$ and $h$ on $(0,1)$
and the last inequality we can conclude the proof.
\cqfd


\subsection{Second kinetic inequality}

A second kinetic formulation not used in this paper is briefly presented.

If $G''>0$, for all weak entropy solutions we have, for all $\phi$ such that $$(u,m) \mapsto
\phi(\ln(u)+g(m/u))$$ is convex, see Section \ref{sh},

\begin{equation*} 
\partial_x \left(\phi(uG(c)) \right)  \leq  0. 
\end{equation*}

Unfortunately these kind of entropy is related with convexity or not of isotherms. Furthermore we
cannot expect to recover all System (\ref{sysad})  (but if $G'' > 0$ we must have the Lax entropy
condition).

In \cite{BGJ2}, we prove that is $w \mapsto p(w)$ is nonnegative and nondecreasing function,
$\phi(w):=p(w)\exp(w)$ is a convex entropy. 

Notice that $w \in \R$ in contrast to $c \in [0,1]$. Let us introduce some notations before
exhibiting a new family of entropies:

\begin{equation*} 
     \zeta_+= \max(\zeta,0)\geq 0,\quad \zeta_-= \min(\zeta,0)\leq 0,
\end{equation*}
 
\begin{eqnarray*}  \label{pw}
 p(w,\zeta)= (w-\zeta)_+ + \zeta_- = 
                  \left \{ \begin{array}{cl} 
                  0  & \mbox{ if } 0 <w< \zeta \mbox{ or }\zeta >0 > w,  \\  
     
                w-\zeta  & \mbox{ if } 0 < \zeta < w, \\
                  \zeta  & \mbox{ if } w < \zeta < 0, \\                        
                w   & \mbox{ if }   \zeta < w< 0  \mbox{ or }\zeta <0 < w, 
             \end{array} \right.
\end{eqnarray*}

\begin{eqnarray*}  \label{kai2}
    \widetilde{\chi}(w,\zeta)= -\partial_\zeta p(w,\zeta)=
                  \left \{ \begin{array}{cl} 
                 1  & \mbox{ if } 0 < \zeta < w, \\
                 -1  & \mbox{ if }  w< \zeta < 0, \\                            
               0  & \mbox{else.}
             \end{array} \right.
\end{eqnarray*}

Furthermore $w(t,x)=\int_\R \widetilde{\chi}(w,\zeta) d\zeta.$
 So, we deduce easily following second kinetic formulation.

\begin{theorem}[Second kinetic formulation]~\\
If $G'' \geq 0$, $(u,c)$ is a weak entropy solution of System (\ref{sysad}), then there exists a
nonnegative measure $\mu(t,x,\zeta)$ on $\R^2_+\times  (0,1)$ such that: 
\begin{equation} \label{eq2kf}
\partial_x \left( uG(c) \, \widetilde{\chi}(w,\zeta) \right)  =\partial_\zeta\mu(t,x,\zeta). 
\end{equation}
Furthermore we have the a priori bound for all $t>0, X>0$:
\begin{equation*}
 \int_0^X\int_\R \mu(t,x,\zeta) \, d\zeta  \, dx  = \int_0^X u \, G(c) \, w(t,x) \, dx. 
\end{equation*}
\end{theorem}

 \begin{remark}
 In the case of one inert gas with ammoniac or water vapor $G''<0$,  (\ref{eq2kf}) is valid but
with a non positive measure $\mu$.  
\end{remark}


\section{A Kinetic scheme for the PSA System}\label{ksPSA}


In this section a kinetic scheme related to the kinetic formulation  \eqref{recip} is proposed . If the velocity $u$ is frozen then the kinetic formulation seems only related to the concentration $c$. Thus a kinetic scheme for scalar conservation laws can be used except that the macroscopic variable $c$ appears in the kinetic velocity.  Then, an important step is to update the velocity. To be consistent with the PDE, the second equation of  PSA System \ref{sysad} is used.
This scheme, presented in Section \ref{ssks}, enjoys some mathematical properties: maximum principle in Section \ref{ssksPM}, $BV$ estimates in Section \ref{ssksBV}  and the scheme satisfies entropy inequality in Section \ref{ssksEntropy}.  Finally the scheme is tested with exact solutions of some Riemann problems and many isotherms in Section
\ref{numex}.

\subsection{The Kinetic scheme} \label{ssks}

Let $T>0$ be the duration of the simulated process. The interval $[0,T]$ is divided in $N$ meshes $]t_{i-1/2},t_{i+1/2}[$ ($t_{1/2}=0$)
with same length $\Delta t=T/N$, centered in $t_i$. At each time step a new spatial mesh $[x_n,x_{n+1}]$ ($x_0=0$) with lenght $\Delta x^n$ is defined, according to some CFL type condition. 
The discrete unknowns $(u_i^n,c_i^n)$ with $n\in\N^*$ and $1\leq i\leq N$ are the
approximations of the velocity and the concentration, respectively, at  $x=x_n$ in the temporal mesh
$Te_i=]t_{i-1/2},t_{i+1/2}[$.\\
The initial and boundary data are taken in account setting:

\begin{equation}\label{ibd}
c_0^n=\int_{x_{n}}^{x_{n+1}}c_0(x)dx,\, \ 0\leq n\leq N-1,\qquad u_i^0=\int_{Te_i}u_b(t)dt,\quad c_i^o=\int_{Te_i}c_b(t)dt,\, i\in\N.
\end{equation} 

Let $I\geq 2$ besome integer. Being given $(u_i^n,c_i^n)$ for $i=1,\cdots,I$,  we denote $u^n$
(resp. $c^n$) the piecewise
function with value $u_i^n$ (resp. $c_i^n$) on $]t_{i-1/2},t_{i+1/2}[$ and we introduce the
following transport equation related to \eqref{recip} on $[x_n,x_{n+1}[$ (with $x$ as the evolution variable):

\begin{eqnarray}\label{kl}
\partial_x(u^n(t)\,\phi(t,x,\xi)) + \partial_t((a(\xi)+h(c^n(t)))\,\phi(t,x,\xi)) & =
& 0\label{eqkin}\\
\hbox{with}\quad t>0,\,\,x_n\leq x<x_{n+1},\,\,\xi\in\R&& \nonumber\\
\nonumber\\ 
\hbox{and}\quad \phi(t,x_n,\xi)=\chi(c^n(t),\xi)\quad t>0,\,\,\xi\in\R.&&
\label{dataeqkin}
\end{eqnarray} 

This equation is solved numerically, using a standard explicit upwind finite volume scheme. More
precisely, we set, for each $\xi\in\R$, $\chi_i^n(\xi)=\chi(c_i^n,\xi)$ and we denote
$v_i^{n+1}(\xi)$ an approximation of the mean value on  $]t_{i-1/2},t_{i+1/2}[$ of the solution
$\phi(t,x_{n+1}^-,\xi)$ of (\ref{eqkin})-(\ref{dataeqkin}). \\
Thus $\phi_i^{n+1}(\xi)$ is given by the following scheme:

\begin{eqnarray*}
\Delta t\,u_i^n\,\left(\phi_i^{n+1}(\xi)- \chi_i^n(\xi)\right)+a(\xi)\,\Delta
x^n\,\left(\phi_{i+1/2}^n(\xi)-\phi_{i-1/2}^n(\xi) \right)\label{vf1}\label{mks}\\ 
+\Delta x^n\,h(c_i^n)\,\left(\chi_i^n(\xi)-\chi_{i-1}^n(\xi) \right) =0,\nonumber
\end{eqnarray*} 

with

\begin{eqnarray*}\label{vf2}
 \phi_{i+1/2}^n(\xi)=\left\lbrace\begin{array}{rcl}
\chi_i^n(\xi) & \hbox{ if } & a(\xi)\geq 0,\\
\chi_{i+1}^n(\xi) & \hbox{ if } & a(\xi)< 0,
\end{array}\right.
\end{eqnarray*} 

i.e. 
$$a(\xi)\,\phi_{i+1/2}^n(\xi)= a^+(\xi)\,\chi_i^n(\xi)-a^-(\xi)\,\chi_{i+1}^n(\xi)$$ 

where $a^+=\max(a,0)$ and $a^-=-\min(a,0)$.\\
Notice that $\phi_i^{n+1}(\xi)$ is no longer a Gibbs equilibrium, i.e. a $\chi$ function. We recover
such an equilibrium setting $c_i^{n+1}=\int_\R \phi_i^{n+1}(\xi)\,d\xi$ and thus getting
$\chi_i^{n+1}$.\\
At the macroscopic level, integrating (\ref{kl}) with respect to the kinetic variable $\xi$ we
get:

\begin{equation}\label{ks1}
c_i^{n+1}=c_i^n - \lambda_i^n\left\lbrace (A^+(c_i^n) - A^+(c_{i-1}^n)) -  (A^-(c_{i+1}^n) -
A^-(c_i^n)) + h(c_i^n)\,(c_i^n - c_{i-1}^n )\right\rbrace 
\end{equation} 

where we have set, as long as  $u_i^n > 0$:
\begin{equation}\label{coeff}
\lambda_i^n=\frac{\Delta x^n}{u_i^n\,\Delta t},\qquad A^\pm(c)=\int_0^c a^\pm(\xi)\,d\xi.
\end{equation}

Finally,  we update the velocity $u$ applying a classical finite difference scheme to the second equation of System \eqref{sysad}:
\begin{equation}\label{ks2}
 u_i^ {n+1}=u_i^n - \frac{\Delta x^n}{\Delta t}\,(h(c_{i+1}^n)-h(c_{i-1}^n)).
\end{equation} 

Notice that(\ref{ibd})-(\ref{coeff})-(\ref{ks1})-(\ref{ks2}) allow to compute
$(u_i^{n+1},c_i^{n+1})$  for $i=1,\cdots,I-1$ only,  because the sign of $a$ is not a priori known, thus the
computation will be effective in all the columns for $0\leq t\leq T=(I-M+1)\,\Delta t$ if $M$ is the
number of spatial meshes. We have thus to extend the temporal domain  and the functions $u_b,\,c_b$ in a suitable way, $M$ being estimated following Remark \ref{versionW}.

In the sequel we note (KS) the kinetic scheme defined by 
(\ref{ibd})-(\ref{coeff})-(\ref{ks1})-(\ref{ks2}).

\subsection{$L^\infty$ estimates}   \label{ssksPM}

\begin{proposition}
 Assume that $u_b\geq\alpha$ for some constant $\alpha>0$. As long as  $u_i^n>0$, if $\Delta x^n$  satisfies at each time step
the CFL  type condition
\begin{equation}\label{cfl}
(\norminf{h}+2\,\norminf{a})\,\Delta x^n\leq \min_{i}(u_i^n)\,\Delta t ,
\end{equation} 
thenwe have the following $L^\infty$ estimates:
\begin{equation}\label{Linfty}
 0\leq \phi_i^n\leq 1,\qquad 0\leq c_i^n\leq \max\{ \norminf{c_b},\,\norminf{c_0}\}\leq 1.
\end{equation} 
\end{proposition}

In (\ref{cfl}), the $L^\infty$ norms are relative to $[0,1]$.

\pr writing (\ref{vf1})  under the form
$$\phi_i^{n+1}(\xi)=\left(1-\lambda_i^n\,\left(\abs{a(\xi)} + h(c_i^n) \right)   \right)\,\chi_i^n(\xi)
+ \lambda_i^n\,(a^+(\xi)+h(c_i^n))\,\chi_{i-1}^n(\xi) + \lambda_i^n\,a^-(\xi)\,\chi_{i+1}^n(\xi) $$

we obtain $\phi_i^{n+1}(\xi)$ as a convex combination of $\chi_{i-1}^n(\xi)$, $\chi_i^n(\xi)$ and
$\chi_{i+1}^n(\xi)$ as soon as (\ref{cfl}) is satisfied. Thus the first inequality of
(\ref{Linfty}) holds by induction. In the same manner we can write (\ref{ks1}) as

$$c_i^{n+1}= \left(1-\lambda_i^n\,\left(p_{i-1/2}^{n,+}+p_{i+1/2}^{n,-} + h(c_i^n) \right)  
\right)\,c_i^n + \lambda_i^n\,(p_{i-1/2}^{n,+} + h(c_i^n))\,c_{i-1}^n +
\lambda_i^n\,p_{i+1/2}^{n,-}\,c_{i+1}^n$$

where $ p_{i-1/2}^{n,+}=\frac{A^+(c_i^n)-A^+(c_{i-1}^n)}{c_i^n-c_{i-1}^n}$ and $
p_{i+1/2}^{n,-}=\frac{A^-(c_{i+1}^n)-A^-(c_i^n)}{c_{i+1}^n-c_i^n}$ satisfy:
$$0\leq  p_{i-1/2}^{n,+}\leq \norminf{a^+} \hbox{ and } 0\leq p_{i+1/2}^{n,-}\leq \norminf{a^-} $$
and the second inequality of (\ref{Linfty}) holds by induction since
$\norminf{a^+}+\norminf{a^-}\leq 2\,\norminf{a}$.\cqfd

\begin{figure}[H]
\centering
\includegraphics[scale=0.6]{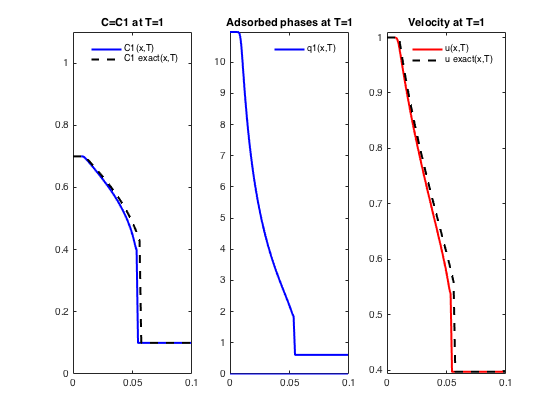}
\caption{adsorption step with the BET isotherm. Exact and computed solutions at $t=1$ with $u$ computed with the Riemann invariant $u\,G(c)$ \label{betW}}
\end{figure}

\begin{remark}\label{versionW}
The previous result is not fully satisfactory because  we are not currently able to give a positive lower bound for $u$ as with the Godunov scheme (\cite{BGJ2}). With the assumption $\inf u_b=\alpha >0$, if we choose to update the velocity using the Riemann invariant $W$, that is setting

\begin{equation}\label{W}
u_i^{n+1}\,G(c_i^{n+1})=u_i^{n}\,G(c_i^{n})
\end{equation}
we get immediately, by induction : $u_i^n\geq \alpha\,\frac{\inf_{[0,1]}G}{\norminf{G}} $. Then we can use the uniform CFL condition:
\begin{equation*}\label{cflu}
\norminf{G}\, (\norminf{h}+2\,\norminf{a})\,\Delta x\leq \alpha\,\Delta t\,\inf_{[0,1]}G ,
\end{equation*} 
and of course \eqref{Linfty} holds. It is not clear a priori that this variant of the kinetic scheme is able to correctly solve the Riemann problem, but the numerical tests show that the behavior is quite satisfactory from this point of view: see Fig. \ref{betW}.
\end{remark}

\subsection{$BV$ estimates}\label{ssksBV}
This scheme is  a TVD scheme.
Let us define the total variation (with respect to the time variable) of $(c_i^n)_{i\in\N}$ by
$$TV^n(c)=\sum_{i=0}^\infty \Big|\Delta c_{i+1/2}^n\Big|,\quad\hbox{ with }\Delta c_{i+1/2}^n=c_{i+1}^n- c_i^n.$$

\begin{proposition}\label{tvd}
The kinetic scheme (KS) is total variation diminishing, that is $$\forall n\in\N\quad TV^{n+1}(c)\leq  TV^n(c) + TV c_b, $$
under the (CFL) condition (\ref{cfl}).
\end{proposition}

\pr we have just to show that $c_i^n$ may be written under an incremental form (see for instance \cite{GR91, GR96}). Now we have:
$$c_i^{n+1}= c_i^n + A_{i+1/2}^n \Delta c_{i+1/2}^n - B_{i+1/2}^n \Delta c_{i-1/2}^n$$
with 
$$A_{i+1/2}^n=\lambda_i^n p_{i+1/2}^{n,-}\quad \hbox{ and } \quad
B_{i+1/2}^n=\lambda_i^n\,\left(h(c_i^n)+ p_{i-1/2}^{n,+} \right).$$
So it is easy to verify that under the condition (\ref{cfl}) we have
$$A_{i+1/2}^n\geq 0,\quad B_{i+1/2}^n \geq 0,\quad \hbox{ and } A_{i+1/2}^n+B_{i+1/2}^n\leq 1$$

which ensures the incremental form for $c_i^n$ and thus concludes the proof.\cqfd

\subsection{Discrete entropy inequalities} \label{ssksEntropy}

An important requirement for  the scheme is to satisfy some entropy inequalities. This is possible with the choice of updating the velocity $u$. 

\begin{proposition}
 The discrete  unknowns $(u_i^n,c_i^n)$ satisfy the following discrete entropy inequality, where $\psi$ is any real smooth function such that $\psi''\geq 0$:
 \begin{equation}\label{dei}
  u_i^{n}\psi(c_i^{n+1})  \leq  u_i^n\psi(c_i^n) - \frac{\Delta x^n}{\Delta t}\,(\phi_{i+1/2}^n-\phi_{i-1/2}^n) -
   \frac{\Delta x^n}{\Delta t}\,h(c_i^n)\,(\psi(c_i^n)-\psi(c_{i-1}^n))
 \end{equation}
 where $\phi_{i+1/2}^n=\phi^+(c_i^n) - \phi^-(c_{i+1}^n)$ with $\phi^\pm=\int\psi'\,a^\pm$.
 
\end{proposition}

\pr  Inequality (\ref{dei}) arises  mutiplying (\ref{mks}) by $\psi'(\xi)$, integrating over
$\R$ with respect to $\xi$ and applying  Brenier's lemma (see Section \ref{kfPSA}) thanks to the
$L^\infty$ estimate (\ref{Linfty}) on $f_i^{n+1}(\xi)$.\cqfd


\section{The kinetic scheme and the Riemann Problem}\label{numex}

In this section our aim is to test the ability of the kinetic scheme to select the entropy solution of the Riemann problem for various choices 
of isotherms following \cite{DCRBT88}. In particular, the BET isotherm with one inflexion point  leads to composite waves which we will see that they are properly calculated by the scheme. 
For self contain we recall the solution of the Riemann Problem (see \cite{BGJ2}), 
moreover  the hyperbolicity condition $u>0$ is ensured for this solution by Proposition \ref{pu>0}.

\subsection{Exact solution of the Riemann Problem}

We consider the following Riemann problem:

\begin{eqnarray}
\left\{   \begin{array}{ccc}  \vspace{2mm}
\partial_x u +\partial_t h(c)&=&0,\\
\partial_x  (u c)+\partial_t I(c)  & = & 0,
\end{array}\right.\label{sysPR}\\
c(0,x)=c^- \in [0,1],  \quad x > 0, &\quad
&
\left\{
\begin{array}{ccl}
c(t,0) &=&c^+   \in [0,1],\\
u(t,0)&=&u^+ > 0,
\end{array}
\right.
t>0 \label{dataPR}
\end{eqnarray}

and we search a selfsimilar solution, i.e. : $ c(t,x)=C(z)$,  $u(t,x)=U(z)$ with $z=\frac{t}{x} >0$.

\begin{figure}[H]
\centering
\includegraphics[scale=0.5]{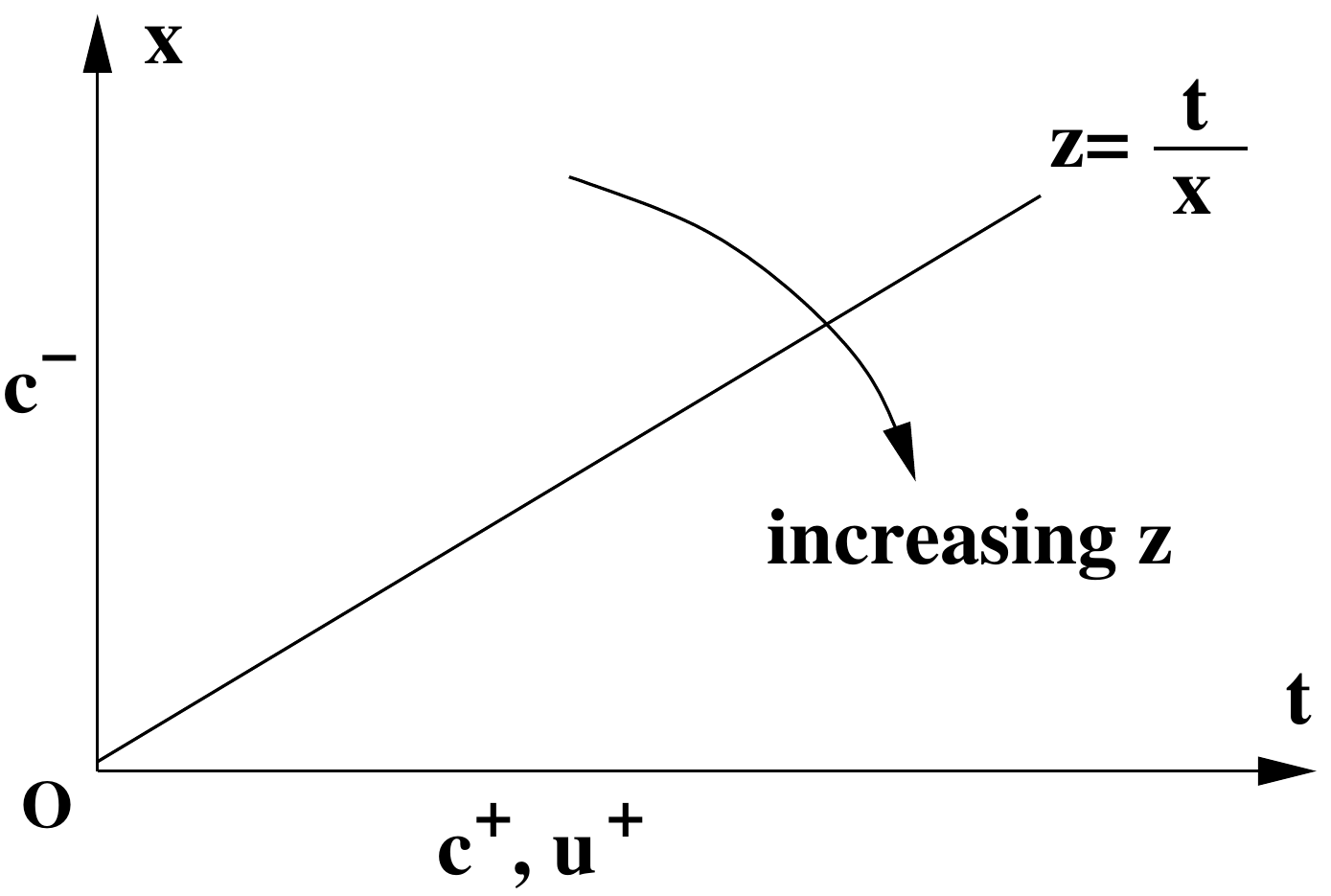}
\caption{data for the Riemann problem \label{drp}}
\end{figure}
The exact solution is computed using the following results stated for instance in \cite{BGJ2}.

\begin{proposition}[Rarefaction waves] \label{rw}~\\
Any smooth non-constant self-similar solution $(C(z), U(z))$ of (\ref{sysPR}) in an open domain\\ $\Omega=\{0\leq\alpha<z<\beta\}$   where $f''(C(z))$ does not vanish, satisfies:

\begin{equation*}\label{Czab}
\frac{d C}{ dz}  = \frac{H(C)}{z\,f''(C) },
\end{equation*}

\begin{equation*}\label{Uzab}
U(z) = \frac{H(C)}{z}.
\end{equation*}

In particular, $ \frac{d C}{ dz}$ has the same sign as $ f''(C)$.
\end{proposition}

\begin{corollary}\label{corrw}
Assume for instance  that $0 \leq a<c^-<c^+<b \leq 1$ and $f''>0$ in $]a,b[$. Then the only smooth self-similar solution of (\ref{sysPR})-(\ref{dataPR}) is such that :

\begin{equation*}\label{Cz}
\left \{\begin{array}{cccr}
C(z)&=&c^- ,&   0  <  z  <  z^- ,\\
\frac{d C}{ dz} & = &
  \frac{H(C)}{z\,f''(C) },&\;   z^-  <   z  < z^+ \\
C(z) & =&   c^+, &   z^+  <  z,
\end{array}  \right.
\end{equation*}

where

$ z^+=\frac{H(c^+)}{u^+} $, $z^-=z^+\,e^{-\phi(c^+)}$ with $\phi(c)=\int_{c^-}^c \frac{f''(\xi)}{H(\xi)}\,d\xi$. Moreover $U$ is given by:

\begin{equation*}\label{Uz}
\left \{  \begin{array}{cccr}
U(z)& = &   u^0 ,&   0  <   z  <   z^- ,\\
U(z)& =&\frac{H(C(z))}{z}, & \;  z^-   <   z   < z^+, \\
U(z)& =&   u^+ &   z^+  <   z.
 \end{array}  \right.
\end{equation*}

where $u^0=\frac{H(c^-)}{z^-}$.
\end{corollary}

\begin{remark}
It appears that $c$ is always monotone along a rarefaction wave but no longer $u$ because the sign of $h'$ may change. Indeed the Riemann invariant $w=\ln u+g(c)$ is constant along such a wave and $g'$, $h'$ have opposite signs. However notice that in the case where one gas is inert,  $u$ is monotone (see also \cite{BGJ1}).
\end{remark}

We are looking now for admissible shocks in the sense of Liu \cite{L76}.  

\begin{proposition}[$\lambda-$shock waves]\label{SW}
If $(c^-,c^+)$ satisfies the following admis\-si\-bility condition equivalent
to the Liu entropy-condition: 

$$\hbox{for all } c \hbox{ between } c^- \hbox{ and } c^+,
\quad \frac{f(c^+)-f(c^-)}{c^+ - c^-}\leq
\frac{f(c)-f(c^-)}{c - c^-},$$
then the Riemann problem (\ref{sysPR})-(\ref{dataPR}) is solved by a shock wave defined as
\begin{equation*}
C(z)=\left\{\begin{array}{ccl}
c^- &\hbox{ if } &  0 < z < s, \\
c^+ & \hbox{ if }&  s<z,
\end{array}\right.
\qquad
U(z)=\left \{\begin{array}{ccl}
u^0 &\hbox{ if }& 0< z < s,\\
u^+ &\hbox{ if } & s < z,
\end{array}  \right.
\end{equation*}
where $u^0$ and the speed $s$ of the shock are obtained through 
\begin{equation} \label{RH}
 u^0([I] - c^-[h])= u^+([I] - c^+[h]), \qquad \qquad s=\dfrac{[h]}{[u]},
\end{equation}
with
$$ [u]=u^+ - u^-, \;\; [h]=h(c^+) - h(c^-), \;\; [I]=I(c^+) - I(c^-).  $$
\end{proposition}

\begin{proposition}\label{discont}
Two states $U^-$ and $U^+$ are connected by a contact discontinuity if and only if $c^-=c^+$ (with of course $u^-\neq u^+$), or $c^-\neq c^+$ and $f$ affine between $c^-$ and $c^+$.
\end{proposition}

Finally, concerning the Riemann problem, we make use of the following wave fan admissibility criterion (see \cite{D00} for instance): the fan is admissible if each one of its shocks, individually, satisfies the Liu shock admissibility criterion.\\
Then, in view of the previous results, we get the solution of the Riemann problem (\ref{sysPR})-(\ref{dataPR}) for  $c$ in a very simple way, similar to the scalar case with flux $f$.

\begin{figure}[H]
\centering
\includegraphics[scale=0.3]{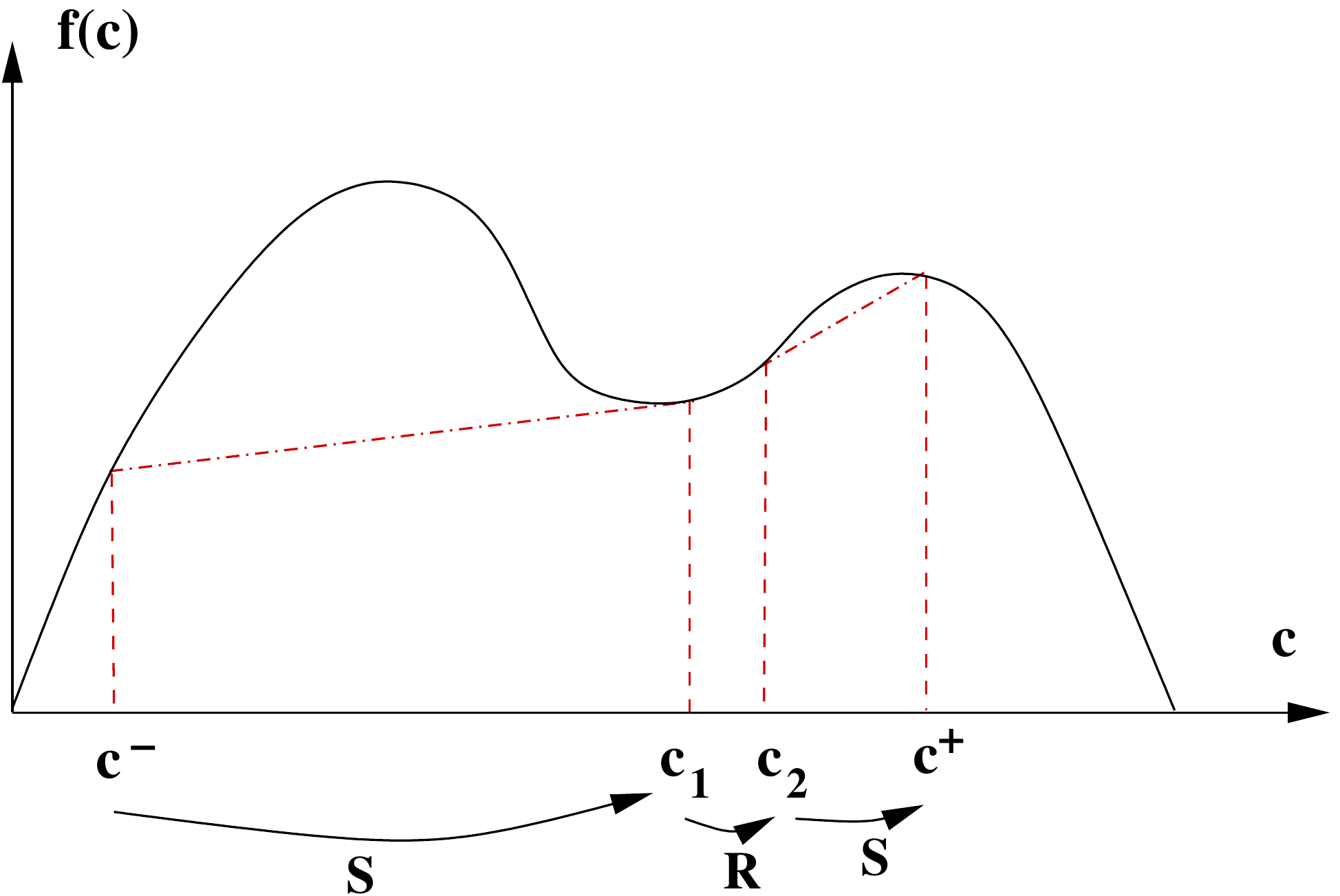}
\includegraphics[scale=0.3]{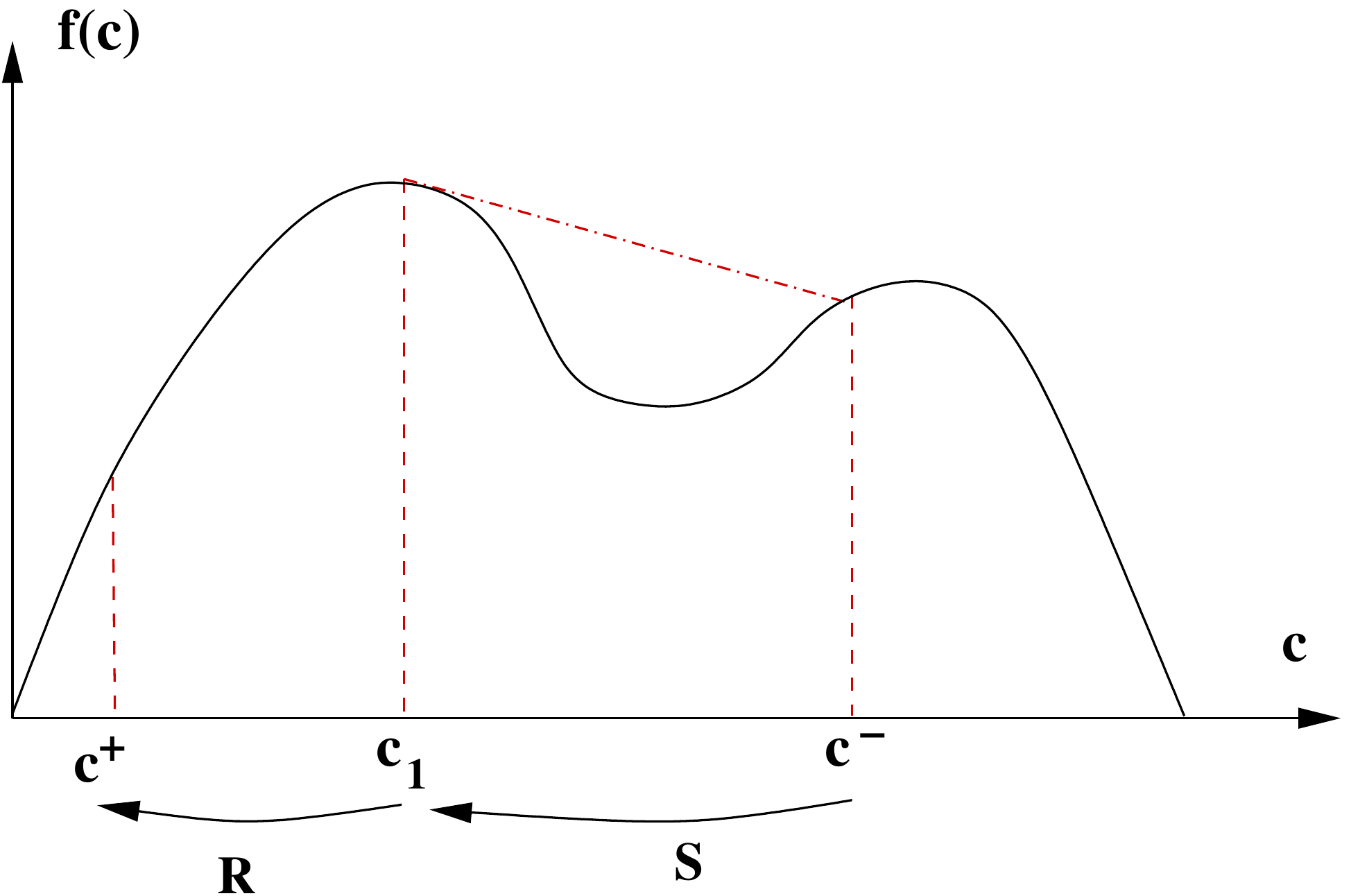}
\caption{ shocks chords are shown as dashed lines. On the left $c^-$ is connected to $c^+$ via a shock (S), a rarefaction wave (R) and a shock. On the right, $c^-$ is connected to $c^+$ via a shock and a rarefaction wave.\label{luc}}
\end{figure}

\begin{description}
\item[Case $c^-<c^+$: ]we consider the lower convex envelope $f_c$ of the function $f$ (see Fig. \ref{luc}, left). On the subintervals where $f$ is strictly convex (then $f=f_c$ ) we get a rarefaction wave according to Corollary \ref{corrw}. Elsewhere we get admissible shock waves (or contact discontinuities if $f$ is affine).
\item[Case $c^->c^+$: ]we use the upper convex envelope $f^c$ (see Fig. \ref{luc}, right) and get rarefaction waves where $f$ is strictly concave.
\end{description}

The subsection is concluded by a new general result about the positivity of the velocity witch improves a similar result in \cite{BGJP}. In other words the region $\{ 0<u,\; 0 \leq c \leq 1\}$ is an invariant domain  for Riemann Problems.  Notice that the positivity of $u$ is mandatory to keep the system hyperbolic, and the velocity can blow up \cite{BGJ5}.

\begin{proposition}\label{pu>0}
Assume that $f$ has a finite number of inflexion points, then the solution of the Riemann Problem with $u^+>0$ involves a positive velocity.
\end{proposition}

\pr with the previous results, the solution of the Riemann Problem consists in a finite sequence of simple waves and it remains to show that the result holds for a simple wave. 
In the case of a  rarefaction wave we have, using \eqref{W},  $u^0=u^+\frac{G(c^+)}{G(c^-)}>0$. 
In the case of a shock wave, we rewrite \eqref{RH} as follows:
$$u^0\,(\underbrace{1+\frac{[q_1]}{[c]}\,(1-c^-)- c^-\,\frac{[q_2]}{[c]}}_{A}) =
u^+\,(\underbrace{1+\frac{[q_1]}{[c]}\,(1-c^+)- c^+\,\frac{[q_2]}{[c]}}_{B}).$$
We have $\frac{[q_2]}{[c]}<0<\frac{[q_1]}{[c]}$, thanks to \eqref{isoth}, and $0\leq c\leq 1$ thus $A,\,B\geq 1$ and $u^0>0$. \cqfd


\subsection{One adsorbable component and inert gas}

In this subsection we compare the exact solution of a Riemann problem with the approximation given by the kinetic scheme in the case of one active gas and one inert gas, with various isotherms. The following numerical examples show the accuracy of the  scheme for  contact discontinuities and composite waves.\\ 
Assume that $q_2^*=0$: the first component is the only active gas. The lenght of the column is $L=0.1$ and $50$ time meshes are used.

\subsubsection{Contact discontinuity}

In this first test case, we use the following isotherm $$q_1^*(c)=K_1\,\frac{c}{1-c}$$
for which the function $f$ is linear: $f(c)=q_1^*(c)\,(1-c)=K_1c$. According to Prop. \eqref{discont}, the Riemann Problem in the $(t,x)$ plane is solved by a contact discontinuity connecting $(c^-,u^-)=(0.2,0.2)$ to $(c^-,u^0)$, with $u^0\simeq 0.10385$, followed by a contact discontinuity (due to the linearity of $f$) connecting  $(c^-,u^0)$ to $(c^+,u^+)=(0.7,0.2)$.\\
In this simulation we have set $K_1=1$. 

\begin{figure}[H]
\centering
\includegraphics[scale=0.6]{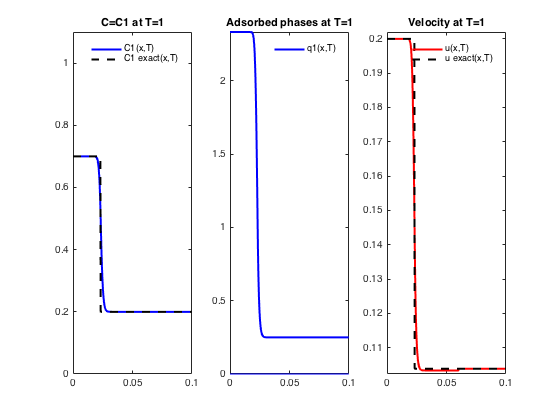}
\caption{ contact discontinuity. Exact and computed solutions at $t=1$ along the column $0\leq x\leq 0.1$ \label{disc}}
\end{figure}

\subsubsection{Adsorption step with the BET isotherm: combined waves}

The so-called BET isothem, in our adimensional variables, is given by:
\begin{equation*}\label{betiso}
q_1^*(c)=\frac{Q\,K\,c}{(1+(K-\frac{1}{c_s})\,c)\,(1-\frac{1}{c_s})}
\end{equation*}
with $0<c_s<1$.\\
In this simulation we have set $Q=1$, $K=10$ and $\frac{1}{c_s}=1.3$. These choices are done in order to obtain a corresponding function $f$ with an inflexion point more easily visible in Fig. \ref{Riemann-BET} below. The Riemann Problem in the $(t,x)$ plane is solved by a contact discontinuity connecting $(c^-,u^-)=(0.1,1)$ to $(c^-,u^0)$, with $u^0\simeq 0.39701$, followed by a shock connecting  $(c^-,u^0)$ to $(c^*,u^*)\simeq (0.41546,0.54985)$ and a rarefaction connecting $(c^*,u^*)$ to $(c^+,u^+)=(0.7,1)$.

\begin{figure}[H]
\centering
\includegraphics[scale=0.1]{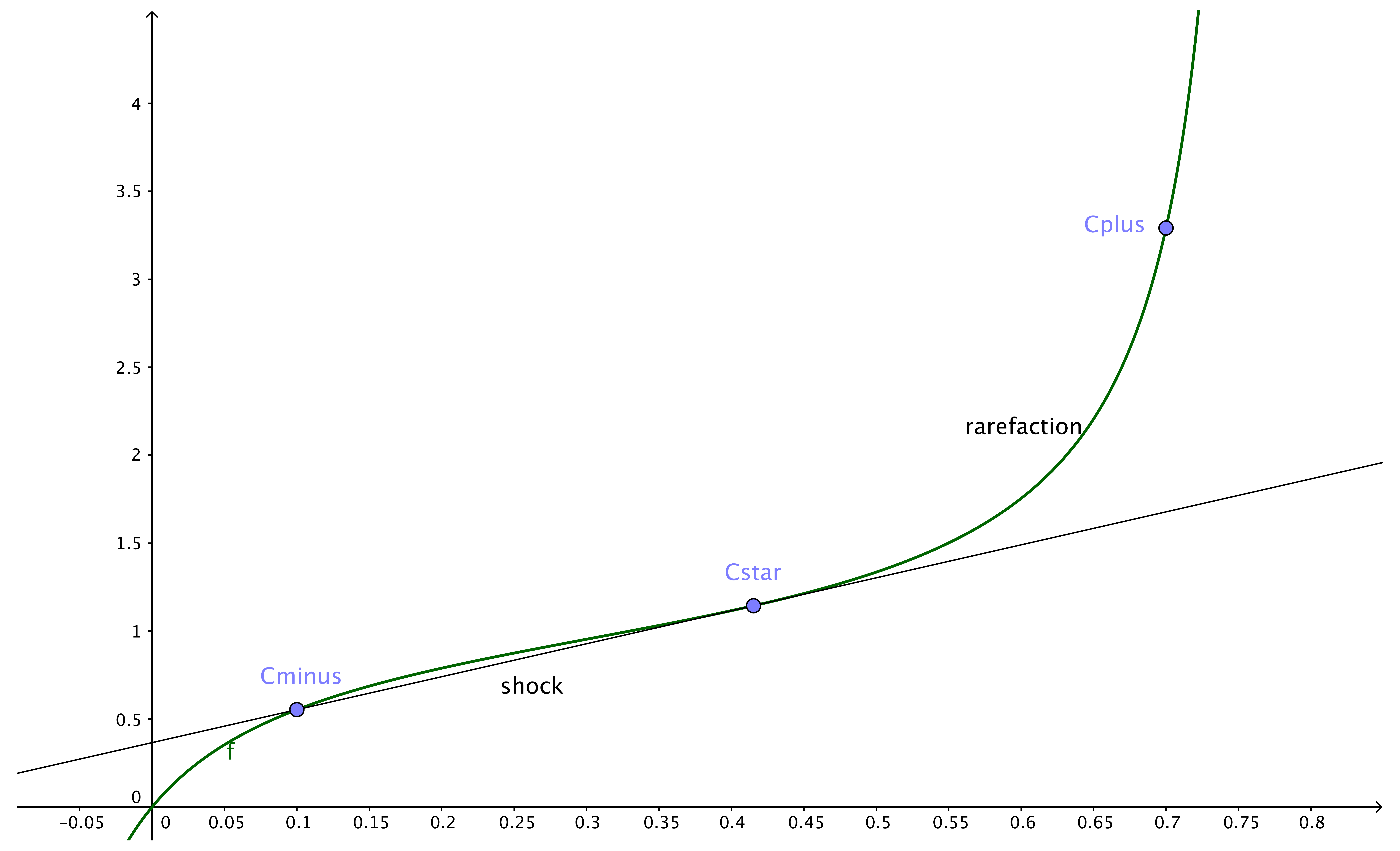}
\caption{solution of the Riemann Problem for $c$ \label{Riemann-BET}}
\end{figure}

\begin{figure}[H]
\centering
\includegraphics[scale=0.6]{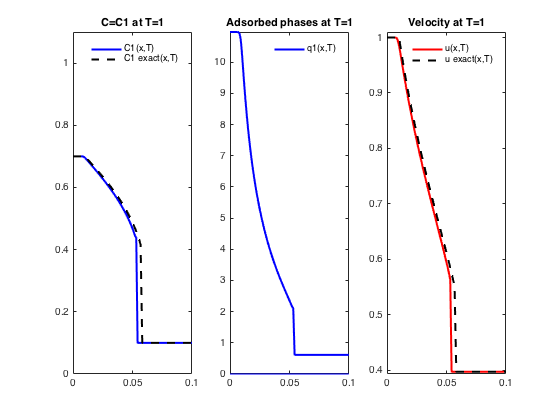}
\caption{adsorption step with the BET isotherm. Exact and computed solutions at $t=1$ along the column $0\leq x\leq 0.1$ \label{bet}}
\end{figure}

We give below the result obtained by updating $u$ through the relation \eqref{W}. It turns out that they are quite similar: the shock and the rarefaction are in both cases correctly computed.

\subsection{Two adsorbable components  with the binary Langmuir isotherm}

In this subsection, we assume that the two gases are active and that the process is driven by the binary Langmuir isotherm:
$$q_i^*(c) = \frac{Q_i K_i c_i}{ 1 + K_1 c_1 + K_2 c_2}, \hbox{ with } K_i >0,\,Q_i>0, \quad i=1,\,2.$$
The following numerical examples show the accuracy of the  scheme for  contact shock and rarefaction waves.\\

The lenght of the column is $L=0.1$ and we used $50$ time meshes. In this simulation,  we have set $Q_1=Q_2=1$, $K_1=10$ and $K_2=30$ : with these values we get a concave function $f$ (see Fig. \ref{bilang} below).

\begin{figure}[H]
\centering
\includegraphics[scale=0.3]{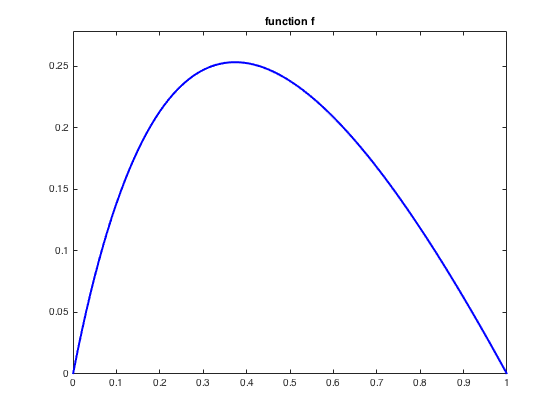}
\caption{function $f$ associated with the binary Langmuir isotherm\label{bilang}}
\end{figure}

The first case, with $c^-<c^+$ (adsorption step) is solved by a contact discontinuity  connecting $(c^-,u^-)=(0.2,0.2)$ to $(c^-,u^0)$, with $u^0\simeq 0.19715$, followed by a shock connecting  $(c^-,u^0)$ to $(c^+,u^+)=(0.7,0.2)$.

\begin{figure}[H]
\centering
\includegraphics[scale=0.6]{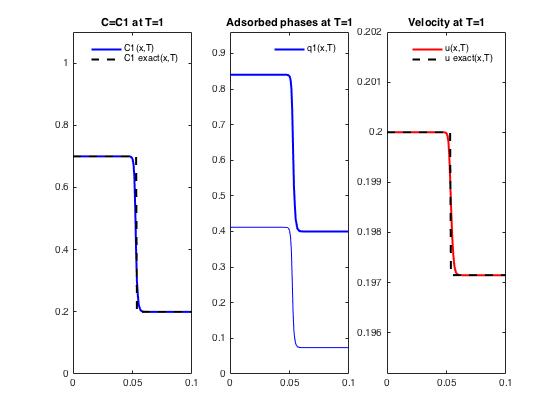}
\caption{adsorption step  with the binary Langmuir isotherm. Exact and computed solutions at $t=1$ along the column $0\leq x\leq 0.1$ \label{bilangadso}}
\end{figure}

The second case, with $c^->c^+$ (desorption step) is solved by a contact discontinuity  connecting $(c^-,u^-)=(0.7,0.2)$ to $(c^-,u^0)$, with $u^0\simeq 0.0.20289$, followed by a rarefaction connecting  $(c^-,u^0)$ to $(c^+,u^+)=(0.2,0.2)$.

\begin{figure}[H]
\centering
\includegraphics[scale=0.6]{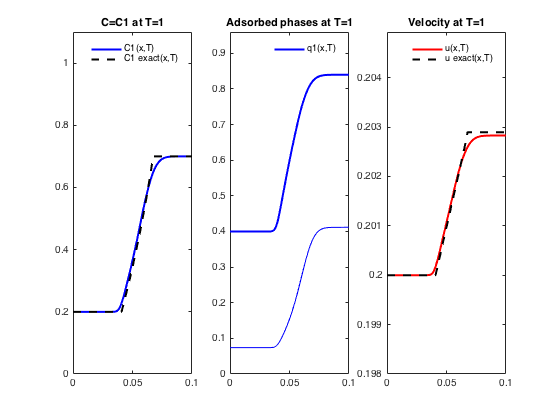}
\caption{desorption step (bottom) with the binary Langmuir isotherm. Exact and computed solutions at $t=1$ along the column $0\leq x\leq 0.1$}
\end{figure}

\section{Conclusion}

We have presented a kinetic formulation of the PSA system, written in an adimensionnal form,  which is used in the context of chemical engineering. This formulation, using an additional real variable,  consits in a single equation which contains, in some sense,  the whole system of two equations and all the entropy inequalities. As a first application, we have  built a kinetic scheme, easy to implement and enjoying good properties (positivity and entropy inequality). It has been tested on the resolution of the Riemann problem in various configurations, including  the case of an isotherm with at least one inflexion point, as the Langmuir isotherm, leading to composite waves. The good agreement with the analytical solution is an argument for convergence and entropic character of the scheme.



\end{document}